\documentclass[10pt]{amsart}

\usepackage{amsmath, amsthm, amscd, amssymb, amsfonts, pstricks}
 
\input amssym.def
\input amssym.tex
\newtheorem{thm}{Theorem}[section]
\newtheorem*{thm*}{Theorem}
\newtheorem{cor}[thm]{Corollary}
\newtheorem*{cor*}{Corollary}
\newtheorem{lem}[thm]{Lemma}

\newtheorem*{con*}{Conjecture}
\newtheorem*{prob*}{Problem}

\theoremstyle{definition}
\newtheorem{defn}[thm]{Definition}

\theoremstyle{remark}
\newtheorem{rem}[thm]{Remark}

\begin{document}
\title[On the (strict) $p$-negative type of finite metric spaces]{Optimal lower bounds
on the maximal $p$-negative type of finite metric spaces}
\author{Anthony Weston}
\email{westona@canisius.edu}
\address{Department of Mathematics {\&} Statistics, Canisius College,
Buffalo, NY 14208}



\begin{abstract}
Determining meaningful lower bounds on the supremal strict $p$-negative type of classes of finite
metric spaces can be a difficult nonlinear problem. In this paper we use an elementary
approach to obtain the following result: Given a finite metric space $(X,d)$ there is a
constant $\zeta > 0$, dependent only on $n = |X|$ and the scaled diameter
$\mathfrak{D} = (\text{diam } X) / \min \{d(x,y) | x \not= y\}$ of $X$ (which we may
assume is $> 1$),
such that $(X,d)$ has $p$-negative type for all $p \in [0,\zeta]$
and strict $p$-negative type for all $p \in [0,\zeta)$. In fact, we obtain
\[
\zeta =  \frac{\ln \bigl( \frac{1}{1 - \Gamma} \bigl)}{\ln \mathfrak{D}}
\text{ where } \Gamma = \frac{1}{2} \Biggl( \frac{1}{\lfloor \frac{n}{2} \rfloor}
+ \frac{1}{\lceil \frac{n}{2} \rceil} \Biggl).
\]
A consideration of basic examples shows that our value of $\zeta$ is optimal provided $\mathfrak{D} \leq 2$.
In other words, for each $\mathfrak{D} \in (1,2]$ and natural number $n \geq 3$, there exists an $n$-point
metric space of scaled diameter $\mathfrak{D}$ whose maximal $p$-negative type is exactly $\zeta$.
It is also the case that the results
of this paper hold more generally for all finite semi-metric spaces since the triangle inequality is
never invoked. Moreover, $\zeta$ is always optimal in the case of finite semi-metric spaces.
\end{abstract}
\maketitle

\section{Introduction}
The notion of $p$-negative type given in Definition \ref{types} (a) was first isolated in the early
1900s through the works of mathematicians such as
Menger \cite{M}, Moore \cite{MOO} and Schoenberg \cite{S1, S2}. One of
the main interests at the time was obtaining characterizations of subsets of Hilbert space up to isometry.
For example, Schoenberg showed that a metric space is isometric to a subset of Hilbert
space if and only if it has $2$-negative type. In the 1960s Bretagnolle \textit{et al}.\ \cite{BDK}
obtained a spectacular generalization of Schoenberg's result to the
category of Banach spaces: A Banach space is linearly isometric to a subspace of some
$L_{p}$-space (for a fixed $p$, $0 < p \leq 2$) if and only if it has $p$-negative type. More recently,
difficult questions concerning $p$-negative type, such as the \textit{Goemans-Linial
Conjecture}, have figured prominently in theoretical Computer Science.
Some monographs and papers which help illustrate the plethora of results and problems in these directions include
Deza and Laurent \cite{DL}, Khot and Vishnoi \cite{KV} (who solved the
\textit{Goemans-Linial Conjecture} negatively), Lee and Naor \cite{LN}, Prassidis and Weston \cite{PW},
and Wells and Williams \cite{WW}.

The related notion of strict $p$-negative type given in Definition \ref{types} (b)
has been rather less well studied than its classical counterpart and most known results deal with the
case $p=1$. Examples of papers which illustrate this concentration on the case $p=1$ include Hjorth
\textit{et al}.\ \cite{HKM, HLM}, Doust and Weston \cite{DW, DW2}, and Prassidis and Weston \cite{PW}.
A theme of this paper is to focus instead on the limiting case $p=0$.
A natural quantification of
the non trivial $p$-negative type inequalities which is given in Definition \ref{NGAP} allows us to
bootstrap the case $p=0$ to obtain the results described in our abstract. This is done
using 
Lagrange multipliers.

In the next section we present some background material
on the notions of (strict) negative type and (strict) generalized roundness. This includes a discussion
of the recently developed notion of the (\textit{normalized}) \textit{$p$-negative type gap}
$\Gamma_{X}^{p}$ of a metric space $(X,d)$. This concept is due to Doust and Weston \cite{DW} and is
formally described in Definition \ref{NGAP}. The informal idea of $\Gamma_{X}^{p}$ is to provide a
natural measure of the degree of strictness of the non trivial $p$-negative type inequalities for any
metric space $(X,d)$ that happens to have $p$-negative type.

The purpose of Section $3$ is to compute the $0$-negative type gap $\Gamma_{X}^{0}$ of each
finite metric space $(X,d)$ exactly. In Theorem \ref{zerogap} we obtain the rather succinct formula
\begin{eqnarray*}
\Gamma_{X}^{0} & = &
\frac{1}{2} \left( \frac{1}{\lfloor \frac{n}{2} \rfloor} + \frac{1}{\lceil \frac{n}{2} \rceil} \right) > 0
\end{eqnarray*}
where $n \geq 2$ denotes $|X|$. Here, given a real number $x$, we are using $\lfloor x \rfloor$ to
denote the largest integer that does not exceed $x$, and $\lceil x \rceil$ to denote the smallest
integer number which is not less than $x$.

The paper culminates in Section $4$ with our main results. We show that every finite semi-metric space
has strict $p$-negative type on an interval of the form $[0,\zeta)$ where $\zeta > 0$
depends only on the cardinality and the scaled diameter of the underlying space.
We also show that $\zeta$ provides an optimal lower bound
on the maximal $p$-negative type of finite metric spaces whose scaled diameter $\mathfrak{D}$ lies in
the interval $(1,2]$. Moreover, in the case of finite semi-metric spaces, we point out that $\zeta$
is \textit{always} optimal. These results appear in Theorem \ref{mainthm}, Remark \ref{mainexample}
and Remark \ref{lastrem}.

Throughout this paper the set of natural numbers $\mathbb{N}$ is taken to consist of all positive integers
and sums indexed over the empty set are always taken to be zero.

\section{Classical (strict) $p$-negative type}
We begin by briefly recalling some theoretical features of (strict) $p$-negative type.
More detailed accounts may be found in the monographs of Deza and Laurent \cite{DL} and
Wells and Williams \cite{WW}. These monographs emphasize the interplay between $p$-negative
type inequalities and isometric embeddings, and (moreover) applications to modern areas such as
combinatorial optimization.

\begin{defn}\label{types} Let $p \geq 0$ and let $(X,d)$ be a metric space. Then:
\begin{enumerate}
\item[(a)] $(X,d)$ has $p$-{\textit{negative type}} if and only if for all natural numbers $k \geq 2$,
all finite subsets $\{x_{1}, \ldots , x_{k} \} \subseteq X$, and all choices of real numbers $\eta_{1},
\ldots, \eta_{k}$ with $\eta_{1} + \cdots + \eta_{k} = 0$, we have:

\begin{eqnarray}\label{ONE}
\sum\limits_{1 \leq i,j \leq k} d(x_{i},x_{j})^{p} \eta_{i} \eta_{j}  & \leq & 0.
\end{eqnarray}

\item[(b)] $(X,d)$ has \textit{strict} $p$-{\textit{negative type}} if and only if it has $p$-negative type
and the associated inequalities (\ref{ONE}) are all strict except in the trivial case
$(\eta_{1}, \ldots, \eta_{k})$ $= (0, \ldots, 0)$.
\end{enumerate}
\end{defn}
\noindent A basic property of $p$-negative type is that it holds on closed intervals of the form $[0,\wp]$.
Namely, if $(X,d)$ is a metric space, then $(X,d)$ has $p$-negative type for all
$p$ such that $0 \leq p \leq \wp$, where $\wp = \max \{ p_{\ast}: (X,d)$ has $p_{\ast}$-negative type$\}$. It is
an open problem if any type of interval property holds for strict $p$-negative type.

It turns out that is possible to reformulate both ordinary and strict $p$-negative type in terms of an
invariant known as \textit{generalized roundness} from the uniform theory of Banach spaces. Generalized roundness
was introduced by Enflo \cite{E} in order to solve (in the negative) \textit{Smirnov's Problem}: Is
every separable metric space uniformly homeomorphic to a subset of Hilbert space? The analogue of this
problem for coarse embeddings was later raised by Gromov \cite{G} and solved negatively by
Dranishnikov \textit{et al}.\ \cite{DGLY}. Prior to introducing generalized roundness in Definition
\ref{GRUNT} (a) we will develop some intermediate technical notions in order to streamline the
exposition throughout the remainder of this paper.

\begin{defn}\label{CUTIE}
Let $q,t$ be arbitrary natural numbers and let $X$ be any set.
\begin{enumerate}
\item[(a)] A $(q,t)$-\textit{simplex} in $X$ is a $(q+t)$-vector $(a_{1}, \ldots , a_{q}, b_{1}, \ldots ,b_{t}) \in X^{q+t}$
whose coordinates consist of $q+t$ distinct vertices $a_{1}, \ldots, a_{q}, b_{1}, \ldots , b_{t} \in X$. Such a
simplex will be denoted $D=[a_{j};b_{i}]_{q,t}$.

\item[(b)] A \textit{load vector} for a $(q,t)$-simplex $D=[a_{j};b_{i}]_{q,t}$ in $X$ is an arbitrary vector
$\vec{\omega} = (m_{1}, \ldots m_{q}, n_{1}, \ldots , n_{t}) \in \mathbb{R}^{q+t}_{+}$ that assigns a positive weight
$m_{j} > 0$ or $n_{i} > 0$ to each vertex $a_{j}$ or $b_{i}$ of $D$, respectively.

\item[(c)] A \textit{loaded} $(q,t)$-\textit{simplex} in $X$ consists of a $(q,t)$-simplex $D=[a_{j};b_{i}]_{q,t}$ in $X$
together with a load vector $\vec{\omega} =(m_{1}, \ldots ,m_{q}, n_{1}, \ldots, n_{t})$ for $D$. Such a loaded simplex
will be denoted $D(\vec{\omega})$ or $[a_{j}(m_{j});b_{i}(n_{i})]_{q,t}$ as the need arises.

\item[(d)] A \textit{normalized} $(q,t)$-\textit{simplex} in $X$ is a loaded $(q,t)$-simplex $D(\vec{\omega})$ in $X$ whose
load vector $\vec{\omega}=(m_{1}, \ldots , m_{q}, n_{1}, \ldots , n_{t})$ satisfies the two normalizations:
\[
m_{1} + \cdots + m_{q} = 1 = n_{1} + \cdots n_{t}.
\]
Such a vector $\vec{\omega}$ will be called a \textit{normalized load vector} for $D$.
\end{enumerate}
\end{defn}
\noindent Rather than give the original definition of generalized roundness $p$ from Enflo \cite{E} we will
present an equivalent reformulation in Definition \ref{GRUNT} (a) that is due to Lennard \textit{et al}.\
\cite{LTW} and Weston \cite{W}. (See also Prassidis and Weston \cite{PW}.)

\begin{defn}\label{GRUNT}
Let $p \geq 0$ and let $(X,d)$ be a metric space. Then:
\begin{enumerate}
\item[(a)] $(X,d)$ has \textit{generalized roundness} $p$ if and only if
for all $q,t \in \mathbb{N}$ and all normalized $(q,t)$-simplices $D(\vec{\omega})
= [a_{j}(m_{j});b_{i}(n_{i})]_{q,t}$
in $X$ we have:
\begin{eqnarray}\label{TWO}
\sum\limits_{1 \leq j_{1} < j_{2} \leq q} m_{j_{1}}m_{j_{2}}d(a_{j_{1}},a_{j_{2}})^{p} +
\sum\limits_{1 \leq i_{1} < i_{2} \leq t} n_{i_{1}}n_{i_{2}}d(b_{i_{1}},b_{i_{2}})^{p} & ~ & \\
\leq \sum\limits_{j,i=1}^{q,t} m_{j}n_{i}d(a_{j},b_{i})^{p}. \hspace*{2.11in} & ~ & ~ \nonumber
\end{eqnarray}

\item[(b)] $(X,d)$ has \textit{strict generalized roundness} $p$ if and only if it has generalized
roundness $p$ and the associated inequalities (\ref{TWO}) are all strict.
\end{enumerate}
\end{defn}
\noindent Two key aspects of generalized roundness for the purposes of this paper are the following equivalences.
Part (a) is due to Lennard \textit{et al}.\ \cite{LTW} and part (b) was later observed by Doust and Weston \cite{DW}.

\begin{thm}\label{REMGR} Let $p \geq 0$ and let $(X,d)$ be a metric space. Then:
\begin{enumerate}
\item[(a)] $(X,d)$ has $p$-negative type if and only if
it has generalized roundness $p$.

\item[(b)] $(X,d)$ has strict $p$-negative type if and only if
it has strict generalized roundness $p$.
\end{enumerate}
\end{thm}

Based on Definition \ref{GRUNT} (a) and Theorem \ref{REMGR} we introduce two numerical
parameters $\gamma_{D}^{p}(\vec{\omega})$
and $\Gamma_{X}^{p}$ that are designed to quantify the \textit{degree of strictness}
of the non trivial $p$-negative type inequalities.

\begin{defn}\label{SGAP} Let $p \geq 0$ and $(X,d)$ be a metric space. Let $q,t$ be natural numbers and
$D=[a_{j};b_{i}]_{q,t}$ be a $(q,t)$-simplex in $X$. Denote by $N_{q,t}$
the set of all normalized load vectors $\vec{\omega}=
(m_{1}, \ldots , m_{q}, n_{1}, \ldots , n_{t}) \subset \mathbb{R}^{q+t}_{+}$ for $D$. Then
the \textit{(normalized)} $p$-\textit{negative type simplex gap} of $D$ is defined to be the
function $\gamma_{D}^{p} : N_{q,t} \rightarrow \mathbb{R}$ where
\begin{eqnarray*}
\gamma_{D}^{p}(\vec{\omega}) & = & \sum\limits_{j,i = 1}^{q,t} m_{j}n_{i}d(a_{j},b_{i})^{p} \\
& ~ & - \sum\limits_{1 \leq j_{1} < j_{2} \leq q} m_{j_{1}}m_{j_{2}}d(a_{j_{1}},a_{j_{2}})^{p}
- \sum\limits_{1 \leq i_{1} < i_{2} \leq t} n_{i_{1}}n_{i_{2}}d(b_{i_{1}},b_{i_{2}})^{p}
\end{eqnarray*}
for each $\vec{\omega}=(m_{1}, \ldots, m_{q}, n_{1}, \ldots, n_{t}) \in N_{q,t}$.
\end{defn}
\noindent Notice that $\gamma_{D}^{p}(\vec{\omega})$ is just taking the difference between the right-hand side
and the left-hand side of the inequality (\ref{TWO}). So, by Theorem \ref{REMGR}, $(X,d)$ has
strict $p$-negative type if and only if $\gamma_{D}^{p}(\vec{\omega}) > 0$ for each normalized
$(q,t)$-simplex $D(\vec{\omega}) \subseteq X$.

\begin{defn}\label{NGAP} Let $p \geq 0$.
Let $(X,d)$ be a metric space with $p$-negative type. We define the
\textit{(normalized)} $p$-\textit{negative type gap} of $(X,d)$
to be the non negative quantity $$\Gamma_{X}^{p} = \inf\limits_{D(\vec{\omega})} \gamma_{D}^{p}(\vec{\omega})$$
where the infimum is taken over all normalized $(q,t)$-simplices $D(\vec{\omega})$ in $X$.
\end{defn}
\noindent Notice (for example) that if $\Gamma_{X}^{p} > 0$, then $(X,d)$ has strict $p$-negative type.
Doust and Weston \cite{DW} have given an example of an infinite metric space to show that the converse
of this statement is not true in general. In other words, there exist infinite metric spaces $(X,d)$ with strict
$p$-negative type and with $\Gamma_{X}^{p} = 0$. It is not at all clear whether the same phenomenon can
occur for finite metric spaces that have strict $p$-negative type.

\begin{rem}
Suppose $(X,d)$ is a metric space with $p$-negative type for some $p \geq 0$. There are two ways in
which we may view the parameter $\Gamma = \Gamma_{X}^{p}$. By definition, $\Gamma$ is the largest non
negative constant so that
\begin{eqnarray}\label{THREE}
\Gamma + \sum\limits_{1 \leq j_{1} < j_{2} \leq q} m_{j_{1}}m_{j_{2}}d(a_{j_{1}},a_{j_{2}})^{p} +
\sum\limits_{1 \leq i_{1} < i_{2} \leq t} n_{i_{1}}n_{i_{2}}d(b_{i_{1}},b_{i_{2}})^{p} & ~ & \\
\leq \sum\limits_{j,i=1}^{q,t} m_{j}n_{i}d(a_{j},b_{i})^{p}. \hspace*{2.11in} & ~ & ~ \nonumber
\end{eqnarray}
for all normalized $(q,t)$-simplices $D(\vec{\omega}) = [a_{j}(m_{j});b_{i}(n_{i})]_{q,t}$ in $X$.
Alternatively, $\Gamma$ is the largest non negative constant so that
\begin{eqnarray}\label{FOUR}
\frac{\Gamma}{2} \left(\sum\limits_{\ell = 1}^{k} |\eta_{\ell}| \right)^{2} +
\sum\limits_{1 \leq i,j \leq k} d(x_{i},x_{j})^{p} \eta_{i} \eta_{j} & \leq & 0.
\end{eqnarray}
for all natural numbers $k \geq 2$,
all finite subsets $\{x_{1}, \ldots , x_{k} \} \subseteq X$, and all choices of real numbers $\eta_{1},
\ldots, \eta_{k}$ with $\eta_{1} + \cdots + \eta_{k} = 0$.
The fact that $\Gamma$ is scaled on the left-hand side of (\ref{FOUR}) simply reflects that the classical
$p$-negative type inequalities are not (by definition) normalized whereas the generalized roundness
inequalities are normalized. See, for example, Doust and Weston \cite{DW}, \cite{DW2}.
\end{rem}

Recall that a \textit{finite metric tree} is a finite connected graph that has no cycles,
endowed with an edge weighted path metric.
Hjorth \textit{et al}.\ \cite{HLM} have shown that
finite metric trees have strict $1$-negative type. Thus it makes sense to try to
compute the $1$-negative type gap of any given finite metric tree. This has been done recently
by Doust and Weston \cite{DW}. However, a modicum of notation is necessary before
stating their result. The set of all edges in a metric tree $(T,d)$,
considered as unordered pairs, will be denoted $E(T)$, and the metric length $d(x,y)$ of any given
edge $e=(x,y) \in E(T)$ will be denoted $|e|$.

\begin{thm}[Doust and Weston \cite{DW}]\label{treegap}
Let $(T,d)$ be a finite metric tree. Then
the (normalized) $1$-negative type gap $\Gamma = \Gamma_{T}^{1}$ of $(T,d)$ is given by the following formula:
\[
\Gamma = \Biggl\{ \sum\limits_{e \in E(T)} |e|^{-1} \Biggl\}^{-1}.
\]
In particular, $\Gamma > 0$.
\end{thm}
\noindent Although strict $1$-negative type has been relatively well studied, properties of strict $p$-negative
type for $p \not= 1$ remain rather obscure and, indeed, there are a large number of intriguing open problems
which beg further investigation. See, for example, Section $6$ of Prassidis and Weston \cite{PW} which lists
some such problems.

\section{Determining the $0$-negative type gap of a finite metric space}
One interesting feature of $0$-negative type is that the metric gets
\textit{forgotten} in the families of inequalities (\ref{ONE}) and (\ref{TWO}) since we will have
$d(x,y)^{0} = 1$ for all $x,y$ with $x \not=y$.
In this section we shall concentrate on the limiting case $p=0$ of negative type in the context of
finite metric spaces.

Forgetting the metric allows the exact computation
of the $0$-negative type gap $\Gamma_{X}^{0}$ for each finite metric space $(X,d)$.
This is done in Theorem \ref{zerogap}.
And, as one would expect, the resulting formula for $\Gamma_{X}^{0} > 0$ depends only on $|X|$.
The most critical computations actually take place in
the following technical lemma. Recall that $N_{q,t}$
denotes the set of vectors $\vec{\omega} = (m_{1}, \ldots, m_{q}, n_{1}, \ldots, n_{t}) \in
\mathbb{R}_{+}^{q+t}$ that satisfy the two constraints $m_{1} + \cdots + m_{q} = 1 =
n_{1} + \cdots + n_{t}$.

\begin{lem}\label{techlem}
Let $(X,d)$ be a metric space and let $q,t \in \mathbb{N}$ such that $q+t \leq |X|$.
Then, for each $(q,t)$-simplex $D=[a_{j};b_{i}]_{q,t} \subseteq X$, we have
\begin{eqnarray*}
\min\limits_{\vec{\omega} \in N_{q,t}} \gamma_{D}^{0}(\vec{\omega})
& = & \frac{1}{2} \left( \frac{1}{q} + \frac{1}{t} \right).
\end{eqnarray*}
In particular, this minimum depends only on $q$ and $t$, and not on the specific vertices
of the simplex $D \subseteq X$.
\end{lem}

\begin{proof}
The overall idea of the proof is to implement Lagrange's multiplier theorem on a large scale.
Consider a given $(q,t)$-simplex $D=[a_{j};b_{i}]_{q,t} \subseteq X$ such that $q+t \leq |X|$.
According to Definition \ref{SGAP}, the $0$-negative type simplex gap $\gamma_{D}^{0}$ is currently
only defined on the constraint surface $N_{q,t} \subset \mathbb{R}_{+}^{q+t}$. Let $\gamma$ denote
the formal extension of $\gamma_{D}^{0}$ to all of the open quadrant $\mathbb{R}_{+}^{q+t}$. In
other words,
\begin{eqnarray*}
\gamma(\vec{\omega}) & = & \sum\limits_{j,i = 1}^{q,t} m_{j}n_{i} \\
& ~ & - \sum\limits_{1 \leq j_{1} < j_{2} \leq q} m_{j_{1}}m_{j_{2}}
- \sum\limits_{1 \leq i_{1} < i_{2} \leq t} n_{i_{1}}n_{i_{2}}
\end{eqnarray*}
for all $\vec{\omega}=(m_{1}, \ldots, m_{q}, n_{1}, \ldots, n_{t}) \in \mathbb{R}_{+}^{q+t}$.
Notice that both $\gamma$ and $\gamma_{|_{N_{q,t}}} = \gamma_{D}^{0}$ depend only on $\vec{\omega}, q$ and $t$,
and not on the specific vertices of the simplex $D \subseteq X$. (This is because we are dealing with the
limiting case of $p$-negative type: $p=0$.)

To complete the proof we introduce two Lagrange multipliers $\lambda_{1}, \lambda_{2}$ and proceed to solve the
system
\begin{eqnarray}\label{FIVE}
~ & ~ &
\left\{ \begin{array}{c}
\frac{\partial}{\partial m_{j}} \biggl( \gamma(\vec{\omega}) -
\lambda_{1} \cdot \sum\limits_{j_{1}=1}^{q} m_{j_{1}} - \lambda_{2} \cdot
\sum\limits_{i_{1}=1}^{t} n_{i_{1}} \biggl) = 0, \,\, 1 \leq j \leq q \\
~ \\
\frac{\partial}{\partial n_{i}} \biggl( \gamma(\vec{\omega}) -
\lambda_{1} \cdot \sum\limits_{j_{1}=1}^{q} m_{j_{1}} - \lambda_{2} \cdot
\sum\limits_{i_{1}=1}^{t} n_{i_{1}} \biggl) = 0, \,\, 1 \leq i \leq t
\end{array} \right.
\end{eqnarray}
subject to the two constraints imposed by the condition $\vec{\omega} \in N_{q,t}$.

The partial derivatives appearing in (\ref{FIVE}) are easily computed and
(together with the two constraints) lead to the following equations:
\begin{eqnarray}
\sum\limits_{i=1}^{t} n_{i} - \sum\limits_{j_{1} \not= j} m_{j_{1}} - \lambda_{1} 
& = & 0 \text{ for all } j, 1 \leq j \leq q \label{SIX} \\
\sum\limits_{j=1}^{q} m_{j} - \sum\limits_{i_{1} \not= i} n_{i_{1}} - \lambda_{2}
& = & 0 \text{ for all } i, 1 \leq i \leq t \label{SEVEN} \\
n_{1} + \cdots + n_{t} & = & 1 \label{EIGHT} \\
m_{1} + \cdots + m_{q} & = & 1 \label{NINE}
\end{eqnarray}
By adding the $q$ equations of (\ref{SIX}) and applying the constraint (\ref{EIGHT}) we obtain
\begin{eqnarray*}
q - (q-1) \sum\limits_{j=1}^{q} m_{j} - q \lambda_{1} & = & 0.
\end{eqnarray*}
Hence $\lambda_{1} = \frac{1}{q}$ by the constraint (\ref{NINE}). So (further) from (\ref{SIX}) and (\ref{NINE})
we therefore have $m_{j} = 1 - \sum\limits_{j_{1} \not= j} m_{j_{1}} = \frac{1}{q} \text{ for all } j, 1 \leq j \leq q$.
Similarly, $n_{i} = \frac{1}{t}$ for all $i$, $1 \leq i \leq t$.

We may now conclude from Lagrange's multiplier theorem that
\begin{eqnarray*}
\min\limits_{\vec{\omega} \in N_{q,t}} \gamma(\vec{\omega})
& = & qt \cdot \frac{1}{qt} - \frac{q(q-1)}{2} \cdot \frac{1}{q^{2}}
- \frac{t(t-1)}{2} \cdot \frac{1}{t^{2}} \\
& = & 1 - \frac{1}{2} \cdot \frac{q-1}{q} - \frac{1}{2} \cdot \frac{t-1}{t} \\
& = & \frac{1}{2} \left( \frac{1}{q} + \frac{1}{t} \right), 
\end{eqnarray*}
thereby establishing the statement of the lemma.

\end{proof}

\begin{thm}\label{zerogap}
Let $(X,d)$ be a finite metric space with cardinality $n = |X| \geq 2$. Then the (normalized)
$0$-negative type gap $\Gamma_{X}^{0}$ of $(X,d)$ is given by the following formula:
\begin{eqnarray*}
\Gamma_{X}^{0} & = &
\frac{1}{2} \left( \frac{1}{\lfloor \frac{n}{2} \rfloor} + \frac{1}{\lceil \frac{n}{2} \rceil} \right).
\end{eqnarray*}
\end{thm}

\begin{proof}
We begin by considering a fixed natural number $m$ such that $2 \leq m \leq n$.
Suppose $D$ is an arbitrary $(q,t)$-simplex in $X$ such that $q+t = m$. We may assume that $q \leq t$
(by relabeling the simplex if necessary). Let $F(q)$ denote $\min\limits_{\vec{\omega} \in N_{q,t}}
\gamma_{D}^{0}(\vec{\omega})$. We proceed to minimize $F$ as a function of $q$.
According to Lemma \ref{techlem}:
\begin{eqnarray*}
F(q)  & =  & \frac{1}{2} \left( \frac{1}{q} + \frac{1}{m-q} \right).
\end{eqnarray*}
Consideration of $F^{\prime}(q)$ shows that $F$ minimizes when $q = \lfloor \frac{m}{2} \rfloor$
(in which case $t = \lceil \frac{m}{2} \rceil$). Finally, as a function of $m$, the expression
\begin{eqnarray*}
(\min F)(m) & = & \frac{1}{2} \left( \frac{1}{\lfloor \frac{m}{2} \rfloor} +
\frac{1}{\lceil \frac{m}{2} \rceil} \right)
\end{eqnarray*}
decreases (strictly) as $m$ increases. We therefore obtain the stated formula for $\Gamma_{X}^{0}$ by
considering the largest possible value of $m$ allowed in this setting: $m=n$.
\end{proof}

\begin{rem}
Consider a finite metric space $(X,d)$ with $n = |X|$. By way of definition, we say that a normalized
$(q,t)$-simplex $D(\vec{\omega}) =[a_{j}(m_{j});b_{i}(n_{i})]_{q,t}$ in $X$ is \textit{extreme for
the zero $0$-negative type gap of $(X,d)$}
if and only if $\gamma_{D}(\vec{\omega}) = \Gamma_{X}^{0}$. The proofs of Lemma \ref{techlem} and
Theorem \ref{zerogap} show that: $D(\vec{\omega})$ is extreme if and only if $q = \lfloor \frac{n}{2} \rfloor$,
$t = \lceil \frac{n}{2} \rceil$, $m_{j} = \frac{1}{q}$ for all $j$ ($1 \leq j \leq q$), and
$n_{i} = \frac{1}{t}$ for all $i$ ($1 \leq i \leq t$).
\end{rem}

The next result makes the point that (unlike all finite metric spaces) no infinite
metric space $(X,d)$ can have a positive $0$-negative type gap $\Gamma_{X}^{0}$.

\begin{cor}
No infinite metric space has a positive $0$-negative type gap.
\end{cor}

\begin{proof}
Let $n \rightarrow \infty$ in Theorem \ref{zerogap} (whence $\Gamma_{X}^{0} \rightarrow 0^{+}$).
\end{proof}

\section{Lower bounds on maximal $p$-negative type of finite metric spaces}
In this section we present the main results of this paper in Theorem \ref{mainthm}
and Remark \ref{mainexample}. In the proof of Theorem \ref{mainthm} we will employ
the notation $(\min F)(m)$ that was introduced in the proof of Theorem \ref{zerogap}.
The overarching strategy is to exploit the \textit{positive} $0$-negative type gap
of each finite metric space $(X,d)$. Recall that the \textit{diameter} of a metric
space $(X,d)$ is the quantity $\text{diam }X = \max\limits_{x,y \in X}d(x,y)$.

\begin{thm}\label{mainthm}
Let $(X,d)$ be a finite metric space with cardinality $n = |X| \geq 3$ (to avoid the trivial case $n=2$).
Let $\mathfrak{D} = (\text{diam } X) / \min \{d(x,y) | x \not= y\}$ denote the scaled diameter of $(X,d)$.
Then $(X,d)$ has $p$-negative type for all $p \in [0, \zeta]$ where
\begin{eqnarray*}
\zeta & = & \frac{\ln \bigl( \frac{1}{1 - \Gamma} \bigl)}{\ln \mathfrak{D}} \text{ with } \Gamma = \Gamma_{X}^{0}.
\end{eqnarray*}
Moreover, $(X,d)$ has strict $p$-negative type for all $p \in [0,\zeta)$.
\end{thm}

\begin{proof}
We may assume that the metric $d$ is not a positive multiple of the discrete metric on $X$. (Otherwise, $(X,d)$ has
strict $p$-negative type for all $p \geq 0$.) Hence $\mathfrak{D} > 1$. We may also assume that
$\min \{d(x,y) | x \not= y\} = 1$ by scaling the metric $d$ in the obvious way, if necessary.
This means that $\mathfrak{D}$ is now the diameter of our rescaled metric space (which we will continue
to denote $(X,d)$). Let $\Gamma$ denote $\Gamma_{X}^{0}$.

Consider an arbitrary normalized $(q,t)$-simplex $D = [a_{j}(m_{j});b_{i}(n_{i})]_{q,t}$ in $X$. Necessarily,
$m = q+t \leq n$. For any given $p \geq 0$, let
\begin{eqnarray*}
L (p) & = & \sum\limits_{j_{1}<j_{2}} m_{j_{1}}m_{j_{2}}d(a_{j_{i}},a_{j_{2}})^{p}
+ \sum\limits_{i_{1}<i_{2}} n_{i_{1}}n_{i_{2}}d(b_{i_{1}},b_{i_{2}})^{p},\,\,{\rm{and}} \\
R (p) & = & \sum\limits_{j,i} m_{j}n_{i}d(a_{j},b_{i})^{p}.
\end{eqnarray*}
By definition of the $0$-negative type gap $\Gamma = \Gamma_{X}^{0}$ we have
\begin{eqnarray}\label{TEN}
L(0) + \Gamma \leq R(0).
\end{eqnarray}
The strategy of the proof is to argue that
\begin{eqnarray}\label{ELEVEN}
L(p) < L(0) + \Gamma
& \text{ and } &
R(0) \leq R(p)
\end{eqnarray}
provided $p > 0$ is sufficiently small. The net effect from (\ref{TEN}) and (\ref{ELEVEN}) is then
$L(p) < R(p)$. Or, put differently, that $(X,d)$ has strict $p$-negative type.
As all non zero distances in $(X,d)$ are at least one we automatically obtain the second inequality
of (\ref{ELEVEN}) for all $p > 0$: $R(0) \leq R(p)$.
Therefore we only need to concentrate on the first inequality of (\ref{ELEVEN}). First of all notice
\begin{eqnarray}\label{TWELVE}
L (p) - L (0) & = & \sum\limits_{j_{1}<j_{2}} m_{j_{1}}m_{j_{2}}
\bigl( d(a_{j_{1}},a_{j_{2}})^{p} - 1 \bigl) \\
& ~ & ~ \nonumber \\
& ~ & + \sum\limits_{i_{1}<i_{2}} n_{i_{1}}n_{i_{2}}
\bigl( d(b_{i_{1}},b_{i_{2}})^{p} - 1 \bigl) \nonumber \\
& \leq & \Biggl(\sum\limits_{j_{1}<j_{2}} m_{j_{1}}m_{j_{2}}
+ \sum\limits_{i_{1}<i_{2}} n_{i_{1}}n_{i_{2}} \Biggl) \cdot \bigl( \mathfrak{D}^{p} - 1 \bigl) \nonumber \\
& \leq & \left( \frac{q(q-1)}{2} \cdot \frac{1}{q^{2}} + \frac{t(t-1)}{2} \cdot \frac{1}{t^{2}} \right)
\cdot \bigl( \mathfrak{D}^{p} - 1 \bigl) \nonumber \\
& = & \Biggl( 1 - \frac{1}{2} \left( \frac{1}{q} + \frac{1}{t} \right) \Biggl)
\cdot \bigl( \mathfrak{D}^{p} - 1 \bigl) \nonumber \\
& \leq & \Biggl( 1 - \frac{1}{2} \left( \frac{1}{\lfloor \frac{m}{2} \rfloor}
+ \frac{1}{\lceil \frac{m}{2} \rceil} \right) \Biggl)
\cdot \bigl( \mathfrak{D}^{p} - 1 \bigl) \nonumber \\
& = & \Bigl( 1 - (\min F)(m) \Bigl) \cdot \bigl( \mathfrak{D}^{p} - 1 \bigl) \nonumber \\
L(p) - L(0) & \leq & \bigl( 1 - \Gamma \bigl) \cdot \bigl( \mathfrak{D}^{p} - 1 \bigl) \nonumber
\end{eqnarray}
by (slight modifications) of the computations in the proofs of Lemma \ref{techlem} and Theorem \ref{zerogap}.
Now observe that
\begin{eqnarray}\label{THIRTEEN}
\bigl( 1 - \Gamma \bigl) \cdot \bigl( \mathfrak{D}^{p} - 1 \bigl) \leq \Gamma
& \text{ iff } & p \leq \frac{\ln \bigl( \frac{1}{1 - \Gamma} \bigl)}{\ln \mathfrak{D}}.
\end{eqnarray}
By combining (\ref{TWELVE}) and (\ref{THIRTEEN}) we obtain the first inequality of (\ref{ELEVEN})
for all $p > 0$ such that $$p < \zeta = \frac{\ln \bigl( \frac{1}{1 - \Gamma} \bigl)}{\ln \mathfrak{D}}.$$
Hence $L(p) < R(p)$ for any such $p$.
It is also clear from (\ref{ELEVEN}), (\ref{TWELVE}) and (\ref{THIRTEEN}) that $L(\zeta) \leq R(\zeta)$.
In total these observations complete the proof of the theorem.
\end{proof}

\begin{rem}\label{mainexample}
The following example illustrates that the value of $\zeta = \zeta(n, \mathfrak{D})$ obtained in Theorem
\ref{mainthm} is sharp provided $\mathfrak{D} \leq 2$. More precisely, for each natural number $n \geq 3$,
there exists an $n$-point metric space $(X,d)$ of scaled diameter $\mathfrak{D} = 2$
whose maximal $p$-negative type is $\zeta(n,2)$.
The idea is to consider a class of simple metrics on certain complete bipartite graphs.
Consider a natural number $n \geq 3$. Set $q = \lfloor n/2 \rfloor$ and $t = \lceil n/2 \rceil$
(whence $q+t = n$).
Let the set $X$ consist of $n$ distinct vertices $a_{1}, \ldots, a_{q}, b_{1}, \ldots, b_{t}$.
Define a metric $d$ on $X$ as follows: $d(a_{j},b_{i})=1$ for all $i$ and $j$,
$d(a_{j_{1}},a_{j_{2}}) = 2$ for all $j_{1} \not= j_{2}$, and $d(b_{i_{1}},b_{i_{2}}) = 2$
for all $i_{1} \not= i_{2}$. This metric space has scaled diameter $\mathfrak{D} = 2$.
Now consider the normalized $(q,t)$-simplex $D=[a_{j}(1/q);b_{i}(1/t)]_{q,t}$ in $X$.
If (\ref{TWO}) holds with exponent $p \geq 0$ for the normalized $(q,t)$-simplex $D$,
then $p$ is subject to the following inequality
\begin{eqnarray*}
\frac{q(q-1)}{2} \cdot \frac{1}{q^{2}} \cdot 2^{p} +
\frac{t(t-1)}{2} \cdot \frac{1}{t^{2}} \cdot 2^{p} 
& = & 2^{p} \cdot \frac{1}{2} \cdot \Biggl( \frac{q-1}{q} + \frac{t-1}{t} \Biggl) \\
& = & 2^{p} \cdot \bigl( 1 - \Gamma \bigl) \\
& \leq & 1
\end{eqnarray*}
where $\Gamma = \frac{1}{2} \left( \frac{1}{q}
+ \frac{1}{t} \right)
= \frac{1}{2} \left( \frac{1}{\lfloor \frac{n}{2} \rfloor}
+ \frac{1}{\lceil \frac{n}{2} \rceil} \right)$.
It follows that the $p$-negative type of $(X,d)$ cannot exceed
\[
\zeta(n,2)  =  \frac{\ln \bigl( \frac{1}{1 - \Gamma} \bigl)}{\ln 2}.
\]
Hence the maximal $p$-negative type of $(X,d)$ is exactly $\zeta(n, 2)$ by Theorem \ref{mainthm}.
We further note that this example can be modified (in the obvious way) to accommodate any scaled diameter
$\mathfrak{D} \in (1,2]$. The obstruction to the case $\mathfrak{D} > 2$ comes from the triangle inequality
(which would be violated in the above setting). However, no such obstruction exists for semi-metric spaces.
Recall that a \textit{semi-metric} is required to satisfy all of the axioms of a metric except (possibly)
the triangle inequality. In this respect we are following Khamsi and Kirk \cite{KK}.
\end{rem}

\begin{rem}\label{lastrem}
In closing we note that Theorem \ref{zerogap} and Theorem \ref{mainthm} hold (more generally) for all
finite semi-metric spaces $(X,d)$. This is because the triangle inequality has played no r\^{o}le in any
of the definitions or computations of this paper. Moreover, in the case of finite semi-metric spaces,
the value of $\zeta(n, \mathfrak{D})$ from Theorem \ref{mainthm} is always optimal. This follows from
the obvious modifications to the examples given in Remark \ref{mainexample}.
\end{rem}

\section*{Acknowledgments}
The author would like to extend special thanks Professor P H Potgieter
(UNISA) for making a large number of incisive mathematical and editorial
comments on preliminary drafts of this paper. The author would also like
to thank Canisius College for the Dean's Summer Research Grant which
made this work possible.

\bibliographystyle{amsalpha}

\end{document}